\documentclass[10pt,twocolumn,twoside]{IEEEtran}
\usepackage{hyperref}
\usepackage{tabu}
\usepackage{graphicx,epsfig,color}
\usepackage{tikz}
\usepackage{pgfplots}
\pgfplotsset{compat=1.14}
\usepackage{csquotes}
\usepackage{amsfonts,amssymb,amsmath,bm}

\usepackage{cite}

\ifCLASSINFOpdf
\else
\fi
\usepackage{subfig}

\newcommand{\norm}[1]{\left\lVert#1\right\rVert}
\begin{document}

\title{Operations- and Uncertainty-Aware Installation of FACTS Devices in a Large Transmission System}


\author{{\bf Vladimir~Frolov} $^{*,a}$ (vladimir.frolov@skolkovotech.ru), {\bf Priyanko~Guha~Thakurta } $^{a}$ (P.GuhaThakurta@skoltech.ru), {\bf Scott~Backhaus } $^{b}$ (backhaus@lanl.gov), {\bf Janusz~Bialek~\IEEEmembership{Fellow,~IEEE} } $^{a}$ (J.Bialek@skoltech.ru)
        and {\bf Michael~Chertkov} \textit{Senior Member, IEEE } $^{a,b,c}$ (chertkov@math.arizona.edu) \\
        $^{a}$ Skolkovo Institute of Science and Technology, Bolshoy Boulevard 30, bld. 1, Moscow, 121205, Russia\\
$^{b}$ Center for Nonlinear Studies and Theoretical Division, T-4, LANL, Los Alamos, NM 87545, USA\\
$^{c}$ Program in Applied Mathematics, 617 N. Santa Rita, University of Arizona, Tucson, AZ 85721
}

\newcommand{\myspecial}[1]{\mathrm{#1}}

\maketitle

\begin{abstract}
Decentralized electricity markets and more integration of renewables demand expansion of the existing transmission infrastructure to accommodate inflected variabilities in power flows. However, such expansion is severely limited in many countries because of political and environmental issues. Furthermore, high renewables integration requires additional reactive power support, which forces the transmission system operators to utilize the existing grid creatively, e.g., take advantage of new technologies, such as flexible alternating current transmission system (FACTS) devices. We formulate, analyze and solve the challenging investment planning problem of installation in an existing large-scale transmission grid multiple FACTS devices of two types (series capacitors and static VAR compensators.) We account for details of AC character of the power flows, probabilistic modeling of multiple-load scenarios, FACTS devices flexibility in terms of their adjustments within the capacity constraints, and long term practical tradeoffs between capital vs operational expenditures (CAPEX vs OPEX).
It is demonstrated that proper installation of the devices allows to do both - extend or improve feasibility domain for the system and also decrease long term power generation cost (make cheaper generation available). Nonlinear, nonconvex, and multiple-scenario-aware optimization is resolved through an efficient heuristic algorithm consisting of a sequence of quadratic programmings solved by CPLEX combined with exact AC PF resolution for each scenario for maintaining feasible operational states during iterations. Efficiency and scalability of the approach is illustrated on the IEEE 30-bus model and the 2736-bus Polish model from Matpower.
\end{abstract}



\section*{Nomenclature}
\begin{IEEEdescription}[\IEEEusemathlabelsep\IEEEsetlabelwidth{$\overline{P}_g$ ($\underline{P}_g$) $\in$ $\mathbb{R}^{N_g}$}]
\item[\underline{Optimization variables:}]
\item
\item[$\overline{\Delta x}$ $\in$ $\mathbb{R}^{N_l}$] Vector of series FACTS capacities.
\item[$\overline{\Delta Q}$ $\in$ $\mathbb{R}^{N_b}$] Vector of shunt FACTS capacities.
\item[$x$ $\in$ $\mathbb{R}^{N_l}$] Vector of line inductances.
\item[$\Delta x$ $\in$ $\mathbb{R}^{N_l}$] Vector of series FACTS settings.
\item[$\Delta Q$ $\in$ $\mathbb{R}^{N_b}$] Vector of shunt FACTS settings.
\item[$V$, $\theta$ $\in$ $\mathbb{R}^{N_b}$] Vectors of buses' voltages and phases.
\item[$P_G$ $(Q_G)$ $\in$ $\mathbb{R}^{N_b}$] Vector of the generators' active (reactive) power injections.
\item [\underline{Parameters:}]
\item
\item[$\overline{P}_G$ ($\underline{P}_G$) $\in$ $\mathbb{R}^{N_b}$] Vector of maximum (minimum) active power generator outputs.
\item[$\overline{Q}_G$ ($\underline{Q}_G$) $\in$ $\mathbb{R}^{N_b}$] Vector of maximum (minimum) reactive power generator outputs.
\item[$\overline{S}$ $\in$ $\mathbb{R}^{2 N_l}$] Vector of the lines' apparent power limits.
\item[$\overline{V}$ ($\underline{V}$) $\in$ $\mathbb{R}^{N_b}$] Vector of the maximum (minimum) allowed voltages.
\item[$P_{D}$ ($Q_{D}$) $\in$ $\mathbb{R}^{N_b}$] Vector of active (reactive) power demands.
\item[$x_0$  $\in$ $\mathbb{R}^{N_l}$] Vector of initial line inductances.
\item[$a$, $K$] Index of the sampled loading scenario and total number of scenarios
\item[$T_{(a)}$] Occurrence probability of scenario $a$.
\item[$N_l$, $N_b$] Number of operational power lines and buses.
\item[$C_{1}$ $\in$ $\mathbb{R}$] Cost per ohm of a series FACTS device.
\item[$C_{2}$ $\in$ $\mathbb{R}$] Cost per MVAr of a static VAR compensator (SVC) FACTS device.
\item[$N_y$ $\in$ $\mathbb{R}$] Planning horizon for the system.
\item[$C(P_G)$] Generation cost function
\item[$M$] Number of segments representing each load duration curve.
\item[$N_i$] Number of scenarios representing each segment
\item[$\underline{\alpha}$ ($\overline{\alpha}$)] Minimum (maximum) loading level.
\item[$\underline{\alpha}_i$ ($\overline{\alpha}_i$)] Minimum (maximum) loading level for segment $i$.
\item[$p_i=w_i$] Occurrence probability of segment $i$.
\item[$\beta$ $(\%)$] Yearly uniform loading growth
\item[$l^0$ $\in$ $\mathbb{R}^{2 N_b}$] Vector of active and reactive loads for the base configuration.

\end{IEEEdescription}
\IEEEpeerreviewmaketitle

\section{Introduction}

\IEEEPARstart{M}{assive} installation of new resources, such as wind and solar, forces power systems to operate in the constrained regimes close to operational limits. When the system becomes constrained, the transmission capacity needs to be expanded in order to provide committed services \cite{ref1}. However, expansion of the transmission network is severely limited by social and environmental constraints, e.g., in Europe. In this manuscript the option of upgrading power transmission by placing and sizing flexible alternating current transmission system (FACTS) devices is discussed. We show that installation of FACTS devices improves and extends system's feasibility domain and reinforces ability of the system to withstand future operational challenges.

FACTS devices are known to be very effective in increasing transmission capacity and improving power system stability. However, their types, locations, and capacities must be allocated properly in order to exploit their benefits. A number of alternative formulations for optimal placement and sizing of FACTS devices have been  proposed. In particular, \cite{ref4,ref8,ref19} were focused on minimizing the operational  cost. An investment cost is minimized in \cite{ref2,ref6}. Reduction of the transmission losses and increase of the power system loadability was analyzed in
\cite{ref3,ref8} and \cite{ref7,ref11,ref18}. Reduction of the load curtailment, and improvement of voltage profile and voltage stability index were discussed in  \cite{ref12,ref6,ref5,ref18}. The resulting mathematical formulations were, typically, made in the existing literature in terms of a Mixed Integer Non-Linear Programming (MINLP). Then MINLP were resolved via the sensitivity analysis \cite{ref12,ref13}; relaxation and/or decomposition to Mixed Integer linear Programming (MILP) \cite{ref16}; and genetic algorithms \cite{ref9,ref10,ref14,ref21}.
The sensitivity-based methods are very efficient for resolving large-scale problems containing only a few indicators aimed at identifying the lines and/or buses that most significantly affect the placement of the FACTS devices. However, such methods cannot optimize device locations, required installed capacities of the devices, and the number of devices required. Genetic algorithms are advantageous for finding global optimal solutions but suffer from extremely slow convergence. Relaxation techniques, utilized to convert MINLP to MILP \cite{ref17}, suffer from the lack of approximation control, whereas decomposition techniques, which consist of substitution of  the original MINLP with a sequence of MILPs, lead to impractically large hierarchies. Main problem with the aforementioned approaches lies in their poor performance in terms of their computational scaling. Methods suggested in the past are handicapped by their ability to resolve problems limited to only relatively small number of nodes (a few hundreds).
Resolution of thousands nodes large models of practical significance with multiple scenarios were not even considered.
Thus a typical approach to resolving the challenge of scaling consisted in relaying on  approximation techniques to simplify modeling of the the line flows \cite{ref6} or substituting AC power flow modeling by DC modeling \cite{ref16}. Unfortunately, these methodologies suffer from the lack of approximation quality assurance thus making them impractical for planning and installation problems of realistic size.
One stand alone computational approach, showing a significant practical promise for resolving realistic size optimization models formulated as MINLP -- Benders decomposition \cite{ref7} and related scenario based decomposition. In \cite{ref7} scalability is demonstrated for the single scenario cases, but for the multiple scenarios it is limited to models with tens of buses. Allowed number of the installed devices should be limited as well.
Finally,  it also worth mentioning that
placement of FACTS devices requires a significant installation cost even for devices with small installed capacities therefore motivating approaches focused on a search for a sparse placement \cite{ref17,ref22,ref23}.


Inspired by the aforementioned prior studies, this manuscript proposes an alternative scalable and AC-based approach to optimal placement and sizing of a sufficiently small number of FACTS devices in a transmission grid.  Highlights of our approach extending and generalizing our prior work based on the DC description \cite{ref22,ref23}  are as follows:
\begin{enumerate}
\item Optimal placement is resolved by incorporating investment and operational variables into the optimization framework simultaneously. Installation of FACTS introduces additional degrees of freedom which are adjusted (along with other operational degrees of freedom) independently for each scenario within the installed capacities.
In other words, placement is resolved by taking into account operational awareness.

\item Capital and operational expenditures are optimized simultaneously. The main advantage of this approach, which to the best of our knowledge has not been discussed in the literature so far,  is that the resulting optimal investment leads to a greater reduction of the operational costs, thus providing additional long-term benefits.

\item Multiple loading scenarios are considered. Scenarios are generated as samples of a probability distribution associated with projected load curves representing seasonal and daily variations. This is in contrast to the existing literature approach, which accounts for a single (usually worst-case) scenario, thus resulting in installation of an expensive device with an unclear role in other cases. Our approach instead finds a single installation (locations and capacities) resolving multiple problems (e.g., overloads, congestion, voltage problems) associated with a multitude of possible scenarios. Also, the optimal settings (within the installed capacities, distinct for different scenarios) are discovered.

\item A novel optimization heuristic algorithm accounting for the full AC model is developed. The algorithm consists of sequential evaluation until convergence (within the preset tolerance) and it includes two substeps at each step. The first substep is an analytic linearization of the basic AC formulas (nonlinear PF equations and nonlinear line constraints) resulting in a quadratic programming (QP) formulation finding investment variables and operational settings for all scenarios. The second substep consists of solving the AC PF for each scenario, thus updating states found in the first substep.

\item The algorithm resolving multiple loading scenarios scales  well; it is capable of finding a solution for large realistic transmission models with thousands of nodes in a computationally acceptable time. The resulting solution produces either an optimal solution or at least a feasible upper bound solution separated from the optimal one by a relatively small gap.

\end{enumerate}

The rest of the manuscript is organized as follows. Optimization framework for finding the optimal location of FACTS devices is set up in Section \ref{sec:optimization}. The algorithm resolving the optimization is explained in Section \ref{sec:algorithm}. The numerical experiments and results are described and discussed in Section \ref{sec:experiment} and Section \ref{sec:planning}. Conclusions and the path forward are presented in Section \ref{sec:conclusions}. Three appendices provide details on the transmission line $\pi$-modeling (Appendix \ref{sec:transmission_modeling}), the design of the operational load scenarios (Appendix \ref{sec:scenarios}), and the procedure for choosing generation configuration to initialize the algorithm (Appendix \ref{sec:init_gen}).

\section{Optimization framework for finding optimal location of FACTS devices}
\label{sec:optimization}

In this section we, first, formulate the problem of optimal placement and sizing of FACTS devices, and then explain and discuss challenges and features of the resulting nonlinear optimization.

Our multi-scenario, operations aware nonlinear and nonconvex optimization problem, allowing installation of the static capacitors at any lines and static var compensators at any loading nodes, is stated, in terms of continuous  variables (no discrete degrees of freedom) for both installation and operational degrees of freedom as follows:
\begin{eqnarray}
\min_{\overline{\triangle x}, \overline{\triangle Q}, y^{(a)}}C_{1}\norm{\overline{\Delta x}}_1+C_{2}\norm{\overline{\Delta Q}}_1 + N_{y}\sum\limits_{a=1}^{K}T_a C(P_G^{(a)}) \label{OPT}
\end{eqnarray}
subject to
\begin{align}
& y^{(a)} = (x,V, \theta, P_G, Q_G,\triangle x,\Delta Q)^{(a)}	&	\forall a	\label{c20}	\\
& x^{(a)}=x_0^{(a)}+\triangle x^{(a)} & \forall a	\label{c2}\\
& P_G^{(a)} = P_{D}^{(a)} + P^{(a)}_{ij}	& \forall a	\label{ac} \\
& Q_G^{(a)} = Q_{D}^{(a)} - \Delta Q^{(a)} + Q^{(a)}_{ij}	& \forall a		\label{reac} \\
& (P_{ij}^{(a)})_i = \sum_{j\sim i}{\Re(S_{ij}^{(a)})} & \forall i, a	\label{acinj}		\\
& (Q_{ij}^{(a)})_i = \sum_{j\sim i}{\Im(S_{ij}^{(a)})} &	\forall i, a	\label{reacinj}		\\
& \underline{P}_G^{(a)}  \leq P_G^{(a)} \leq \overline{P}_G^{(a)}	&	\forall a \label{c8} \\
& \underline{Q}_G^{(a)}  \leq Q_G^{(a)} \leq \overline{Q}_G^{(a)}	&	\forall a \label{cc8} \\
& -\overline{\triangle x} \leq \triangle x^{(a)} \leq \overline{\triangle x} &	\forall a \label{c4}\\
& -\overline{\triangle Q} \leq \triangle Q^{(a)} \leq \overline{\triangle Q} & \forall a \label{c5} \\
& \underline{V}^{(a)} \leq V^{(a)} \leq \overline{V}^{(a)}	&	\forall a \label{c6}		\\
& [\Re(S)^{(a)}]^{T}[\Re(S)^{(a)}]+[\Im(S)^{(a)}]^{T}[\Im(S)^{(a)}]& \nonumber\\
& \leq (\overline{S}^{(a)})^2 & \forall a \label{c9}	
\end{align}
Variables $\overline{\triangle x}, \overline{\triangle Q}$ represent capacities of the newly installed series capacitors (SCs) and static VAR compensators (SVCs). These are investment, i.e. 1st stage, decision variables. $y^{(a)} = (x,V, \theta, P, Q, \triangle x, \triangle Q)^{(a)}$ represent the second-stage decision variables associated with operations under particular scenario (thus labeled by $^{(a)}$). Notice that the second-stage decision variables include scenario-specific settings of the installed FACTS devices. Operational settings (2nd stage decision variables) are constrained by the available capacities (1st stage decision variables).

The objective function in Eq. \eqref{OPT} consists of three terms. The first two terms express the capital investment costs of the installation of the two types of FACTS devices. Guided by the key message from the field of compressed sensing \cite{06Don}, the $l_1$ norm representation is chosen for the investment terms to promote sparsity of the FACTS device placement. The third term in the objective stands for the operational cost introducing explicit dependence on the settings of all the considered scenarios. Here the summation is over $K$ scenarios  accounting for occurrence probabilities of the scenarios ($T_a$) multiplied by the number of service years considered for planning horizon. By including multiple scenarios over the multi-year time horizon, we thus consider operational-aware planning.

It is important to emphasize the main point, and also difficulty, in solving the multi-scenario  problem -- the problem cannot be split into independent optimization problems, each associated with an individual scenario. The scenarios are coupled through the expected operational cost in the objective function. Obviously there always exists at least one scenario for which scenario-dependent decision variables coincide with the respective capacity-related decision variables (otherwise optimal capacities would allow reduction to decrease the investment cost), however for the most part (majority of the scenarios) operational settings are strictly smaller than the capacity.

The constraints supplementing Eq.~(\ref{OPT}) are as follows.
\eqref{c20} describes the state of the system for a given operational scenario. \eqref{c2} bounds actual line inductances, which are adjusted according to the operational value (per scenario) of the installed SC devices. The operational limits are set according to the respective installed capacities, represented by \eqref{c4}. \eqref{ac} and \eqref{reac} state active and reactive power balance at every bus of the network. Elements of vectors $P_G$ ($Q_G$) and $P_{D}$ ($Q_{D}$) are zeros for buses that contain neither generators nor loads. \eqref{acinj} and \eqref{reacinj} represent the net active ($P_{ij} \in \mathbb{R}^{N_b}$) and reactive ($Q_{ij} \in \mathbb{R}^{N_b}$) power injections at the system buses. The term $\Delta Q$ stands for the scenario-dependent SVC shunt compensation constrained by the respective installed capacities in accordance with \eqref{c5}. Limits on active and reactive power generation are expressed by \eqref{c8} and \eqref{cc8}. \eqref{c6} and \eqref{c9} define the voltage and thermal line flow constraints. Here $S_{ij}^{(a)}=S_f^{(a)}$ if $i$ is \enquote{from} end of a line and $S_{ij}^{(a)}=S_t^{(a)}$ if $i$ is \enquote{to} end of a line. See Appendix \ref{sec:transmission_modeling} for details and nomenclature of the $\pi$-line modeling used in this manuscript.

Notice that nonlinearity and nonconvexity of the constraints \eqref{acinj}, \eqref{reacinj}, and \eqref{c9} constitute the main challenge for solving the optimization efficiently. Available nonlinear solvers, such as IPOPT \cite{IPOPT}, scale poorly (exponentially) with increase in the problem size, making the tool useless for optimization over realistic large transmission networks with thousands of nodes.

To complete the optimization problem formulation, one needs to describe how the representative load scenarios are defined. Each of $K$ scenarios, indexed by $a$ in Eqs.~(\ref{OPT}--\ref{c9}), should characterize different loading configurations with occurrence probability. The scenarios may include sampled (typical) configurations and/or contingency (rare) configurations representing different loading regimes. In principle, choosing scenarios appropriate for the optimization (\ref{OPT}) of a grid model is a stand-alone task.  In this manuscript we choose to generate scenarios from the so-called load duration (LD) curve \cite{ref24}. The scenario generation procedure is explained in Appendix \ref{sec:scenarios}.

\section{The Algorithm}
\label{sec:algorithm}

This section describes the algorithm which we suggest to resolve the optimization problem stated in the preceding section. The algorithm consists of the following principal steps:
\begin{enumerate}
\item Scenarios are generated. (In this manuscript the scenario generation scheme based on the Load Duration curve concept is used. See Appendix \ref{sec:scenarios} for details.)

\item Generation configuration is initialized (for each scenario) according to the scheme explained in Appendix \ref{sec:init_gen}.

\item If some of the constraints \eqref{c2}--\eqref{c9} are violated, the operational point of the system is outside of the feasible domain defined by them. The nonlinear constraints \eqref{acinj}, \eqref{reacinj}, and \eqref{c9} are linearized around the current operational point for each scenario. This allows for the construction of a current linearized version of the nonlinear optimization problem \eqref{OPT}-\eqref{c9}. The problem maintain all constraints for all considered scenarios.

\item The resulting linearized problem is solved by QP using the interior point algorithm of the CPLEX solver \cite{CPLEX}.

\item Exact AC PF is solved to update the states obtained in the previous step. This step is needed to prepare a feasible solution for the next iteration.

\item Steps 2--5 are repeated until no constraints remain violated, the target precision is reached, or the maximum allowed number of iterations is reached.
\end{enumerate}
It is important to emphasize that, by design, the algorithm maintains at each iteration a feasible physical state. Also, the algorithm is a heuristic converging to a local minimum which may or may not be the global minimum. An empirical improvement, in terms of convergence to the optimal (or at least reasonable/good) solution, is achieved through experiments with the algorithm's starting point. It was found that initiating the algorithm with the solution corresponding to optimal dispatch ignoring line constraints (see  Appendix \ref{sec:init_gen}) returns satisfactory results. Note that getting not optimal but reasonable solution resolving all the constraints would normally be acceptable (in practice).


The details of the main steps of the algorithm are presented below.

\subsection{Linearization}
\label{subsec:linearize}

Each scenario acts as an input to this part of the optimization heuristics. The operational state for each scenario can be represented as
\begin{equation}
y^{(a)} = (x,V,\theta,P,Q)^{(a)}
\end{equation}
For each scenario, Eq.~\eqref{c9} is linearized using first order Taylor expansion around the current operational point $y_0$:
\begin{equation}
F^{(a)}(y_0^{(a)}) + \nabla F^{(a)}(y_0^{(a)}) (y^{(a)} - y_{0}^{(a)}) \leq (\overline{S}^{(a)})^2
\end{equation}
Where $F(y)$ is a function defining squared absolute value of the apparent power at an end of the line. Similarly, Eqs. \eqref{ac}-\eqref{reacinj} for each scenario can be linearized as
\begin{align}
&\nabla (P-P_{ij}(x,V,\theta,P,Q))^{(a)}(y_0^{(a)}) (y^{(a)}-y_{0}^{(a)}) = 0	\\
&\nabla (Q-Q_{ij}(x,V,\theta,P,Q))^{(a)}(y_0^{(a)}) (y^{(a)}-y_{0}^{(a)}) = 0
\end{align}

Operational variables are adjusted independently for each scenario.\footnote{Operational values of the installed devices are assumed to be flexible (different for each considered scenario).} However, the capacity limits of the devices stay the same (common) for all of the scenarios.

The fact that all controllable parameters of the system stay adjustable/flexible results in degeneracy of the linearized problem. To resolve possible degeneracy, we take advantage of the flexibility associated with redistributing controllable voltages, active powers and reactive powers. Specifically, we introduce the following soft controllable constraint for reactive power dispatch at each QP step:
\begin{equation}
|Q_G^{(a)} - Q_{G_{0}}^{(a)} | \leq \epsilon
\end{equation}
The constraint is accounted for additionally to Eqs.~\eqref{OPT}--\eqref{c9}.

\subsection{Solving the QP problem}
\label{subsec:QP}

The standard QP solver of CPLEX is utilized to solve the linearized problem for all the considered scenarios together. Outputs of this step include values of operational variables for each scenario along with investment variables $\overline{\triangle x}$ and $\overline{\triangle Q}$.

\subsection{Resolving AC PF}
\label{subsec:ACPF}

Notice that solution of the step \ref{subsec:QP} may actually violate the AC-PF balance. Hence the exact AC-PF step is added to maintain a valid/feasible power AC-PF solution at each step of the algorithm.

All together (i.e., in combination) the steps described above provide a feasible solution and resolve a system's contingencies simultaneously and gracefully.

\section{Numerical Experiments: methodology justification }
\label{sec:experiment}

The developed approach is illustrated on examples of the IEEE 30-bus and the 2736-bus Polish models, both available through the Matpower \cite{matpower} software package. The simulations are performed on a Core i7 2600K@4GHz PC with 24 GB of system memory. Both Matlab and Julia implementations, which are comparable in performance, are used. Operational cost is determined using generation cost functions from the given model. Investment cost is a value calculated according to the installation model (capacity $\times$ installation cost). The actual installation cost of a given FACTS device configuration remains a subject for future research.

\subsection{IEEE 30-bus model}
\label{subsec:30experiment}

In this section the advantage of including operational variables and optimizing the expected value of operational cost as an objective along with the investment cost is emphasized. We also illustrate optimality and scalability of the developed heuristics by comparing our algorithm performance with performance of the IPOPT when it is used as a state-of-the-art brute-force solver applied to the exact problem Eq.(\ref{OPT}).

\paragraph{Necessity of including operational variables in determining FACTS device placement}
Let us, first of all, clarify the algorithm's principal advantage of optimizing operational variables simultaneously with the investment variables. The following simulation is performed to demonstrate the actual benefit of the combined use of the capacity variables and the scenario variables.
The base-case system load is increased uniformly by 5\%, which leads to optimal power flow (OPF) infeasibility. Operational cost is not optimized for now ($N_y=0$ in Eq.~\eqref{OPT}). The initial state of the generation (needed to initialize the algorithm) is defined by the first method described in Appendix \ref{sec:init_gen}. Then the two solutions are compared.  One is the actual solution of the developed algorithm with all degrees of freedom available for the optimization. The second solution is constrained by  the same (fixed) generation dispatch (the initial value).
\begin{table}[h!]
\caption{Monetary advantage of considering operational variables, illustrated on the IEEE 30-bus model.}
\centering
\begin{tabu} to \columnwidth { | X[c] | X[c] | X[c] | X[c] | X[c] |}
\hline
 Oper. variables & SVC cap. (MVAr) & SC cap. (\% of init x) & Invest. cost (\$) & Oper. cost (\$/hour)\\
\hline
Fixed & 6.936 (3 SVCs) & 38 (1 SC) & 415930.29 & 614.05  \\
\hline
Free & 1.112 (1 SVC)& 0 (0 SCs) & 55765.71  & 698.24 \\
\hline
\end{tabu}
\label{tab3}
\end{table}

Table \ref{tab3} details the comparison. The significance of accounting for additional available degrees of freedom is obvious. We find out that although the algorithm is able to find a feasible solution in  both cases, the investment cost (objective function value) is  7.5 times smaller in the adjustable generation dispatch case. Emergence of an expensive solution with small operational cost is reported in Table \ref{tab3}. This is an indication that further experiments with tradeoffs between investment and long-term operational costs (fixed in the use cases studied) are to be explored in future studies.

\paragraph{Importance of including operational cost in the objective function}

Another important point to illustrate on the example of the IEEE 30-bus model is that combining the operational cost and the investment cost in the optimization objective is a way to make the optimization relevant to practical planning. Indeed, keeping only the investment cost produces operationally expensive solutions, whereas keeping only the operational cost results in an expensive (and not sparse) installation. Combining the effects of operations and installations in one objective allows for an efficient balance between the two.

To the best of our knowledge, only the investment cost is accounted for in the available literature devoted to placement of FACTS devices.  To mimic this standard approach (accounting for only a single worst-case operational scenario), the optimization horizon is set to zero, $N_y=0$, in Eq.~\eqref{OPT}. Then $N_y$ is increased to 10 years to take effect of operations into account.

\begin{table}[h!]
\caption{Monetary advantages of adding the operational cost to the optimization objective for the IEEE 30-bus model.}
\begin{tabu} to \columnwidth { | X[c] | X[c] | X[c] | X[c] | X[c] | }
\hline
 Plan. horizon ($N_y$) &  Invest. cost (\$) & Oper. cost (\$/hour) &  Total cost (10 years, M\$) & Difference (\%) \\
\hline
0 & 55618.76 & 698.24 & 61.722 & 14.6  \\
\hline
1 &   121838.27 & 616.25 & 55.202 & 2.49  \\
\hline
10 &   249245.37  & 611.98 & 53.858 & 0.0 \\
\hline
\end{tabu}
\label{tab4}
\end{table}

Table \ref{tab4} illustrates the results. Two solutions found for $N_y=1$ and $N_y=10$ reflects sensitivity to the operational cost. Single extreme load configuration is considered (correspondent to a 5\% increase of the load in the base case).  When planning horizon is extended it becomes profitable to invest more into reduction of the operational cost. It is observed that an additional small investment of 200k$\$$ leads to a savings of 14.6$\%$ of the total cost in 10 years. Based on this example, we conclude that accounting for the operational cost in the planning problem is significant. Properly installed FACTS devices allow not only to resolve infeasibility of the loading configuration but also to reduce the generation cost, thus producing lasting long term benefits. It is important to emphasize that by accounting for multiple representable scenarios (as opposed to working with a single scenario) we achieve a much more realistic description of the whole operational space.

\paragraph{Optimality}

To verify performance of the developed heuristics, a single scenario (base case overloaded by 5$\%$) is considered and the results are compared with the ``exact" (Eq.(\ref{OPT})) ones produced by IPOPT (standard, brute-force, nonlinear solver that is still able to handle the 30-bus investment model).

We choose to work with IPOPT because it shows computational advantage over other nonlinear/nonconvex computational platforms applied to problems with structure similar to the one discussed in the paper. As it is shown in \cite{carleton}, IPOPT is on pier or it outperforms Matpower in solving classic ACOPF. An additional advantage of using IPOPT is in its availability within the JUMP/Julia computational environment we rely on. For the actual solver on step \ref{subsec:QP} of our algorithm the CPLEX solver is called (the solver is known to be advantageous for problems with linear constraints).

The optimization horizon is set at 1 year. Table \ref{tab2} shows comparison of heuristics with the benchmark IPOPT. It is observed that (as expected) developed heuristics produce a very tight upper bound for the exact solution, with values of the objective function and structure of the solution that are very close to the exact values.

\begin{table}[h!]
\caption{Comparison of the proposed heuristics with the brute-force IPOPT solution of the exact problem Eq.(\ref{OPT}) for the IEEE 30-bus model.}
\centering
\begin{tabu} to \columnwidth { | X[c] | X[c] | X[c] | X[c] | X[c] |}
\hline
Solver & Bus number & Calculated cap. (MVAr) & Investment cost (k\$) & Total cost (k\$) \\
\hline
IPOPT & 8 & 2.436 & 121.80 & 5520.094  \\
\hline
proposed algorithm & 8 & 2.437 & 121.84  & 5520.159 \\
\hline
\end{tabu}
\label{tab2}
\end{table}

\paragraph{Scalability}

\begin{figure}[h!]
\begin{center}
\resizebox{0.8\columnwidth}{!}{
%
%
\definecolor{mycolor1}{rgb}{0.00000,0.44700,0.74100}%
\definecolor{mycolor2}{rgb}{0.85000,0.32500,0.09800}%
\begin{tikzpicture}

\begin{axis}[%
font=\Large,
width=5.134in,
height=3.329in,
at={(0.861in,0.457in)},
scale only axis,
xmin=0,
xmax=1000,
xlabel style={font=\color{white!15!black}},
xlabel={\huge Number of scenarios},
ymin=0,
ymax=30,
ylabel style={font=\color{white!15!black}},
ylabel={\huge Time (sec*100)},
axis background/.style={fill=white},
legend style={legend cell align=left, align=left, fill=none, draw=none}
]
\addplot [color=mycolor1, dashed, line width=1.0pt, mark=o, mark options={solid, mycolor1}]
  table[row sep=crcr]{%
1	0.06014\\
5	0.15949\\
10	0.33452\\
20	0.92001\\
40	0.71248\\
60	1.03356\\
80	1.40673\\
160	2.75939\\
320	6.11348\\
1000	18.993\\
};
\addlegendentry{\huge Proposed heuristics}

\addplot [color=mycolor2, dotted, line width=4.0pt, mark=square, mark options={solid, mycolor2}]
  table[row sep=crcr]{%
1	0.00175\\
5	0.00996\\
10	0.02132\\
20	0.23366\\
40	1.61591\\
60	5.57057\\
80	28.79944\\
};
\addlegendentry{\huge IPOPT}

\end{axis}
\end{tikzpicture}%
}
\caption{Computational time comparison for 30-bus model.}
\label{timee_30}
\end{center}
\end{figure}
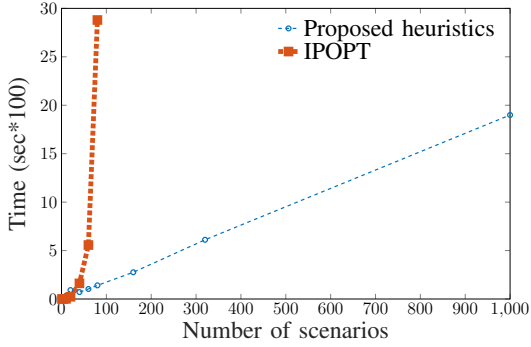

To study how the algorithm scales with the number of scenarios, we pick the base case, increase all loads by $5\%$, and generate $K$ scenarios through the Gaussian sampling procedure associated wtih Eqs.~(\ref{Gauss}--\ref{Gauss-p}), where the rescaled base-case load is set to $l_0$. Scaling analysis of the developed algorithm is illustrated in Fig.~\ref{timee_30}.

It is observed that algorithm handles an increasing number of scenarios very efficiently, solving formulations with a large number of scenarios in time that increases linearly with the number of scenarios (blue line).  Performance of the brute-force IPOPT solver applied to the exact problem in Eq.~\ref{OPT} (shown as thick red line in Fig.~\ref{timee_30}) is limited.

\subsection{2736-bus Polish model}

In this section we extend analysis of the developed algorithm to the case of the Polish grid, which is a practical-size transmission model available as a part of the Matpower package \cite{matpower}. Similar approach as the one tested above on the 30-bus model is followed here.

\paragraph{Necessity of including operational variables in determining FACTS device placement}

The experimental setting stays the same as in the case of the 30-bus model. A single scenario corresponding to the normal operational state (correspondent to the base-case example from  \cite{matpower}) with all the loads homogeneously increased by $5\%$ is considered for resolution. Generation dispatch is defined by the second procedure from Appendix \ref{sec:init_gen}. Table \ref{tab5} illustrates the results. (See Section \ref{subsec:30experiment} for detailed discussion of the experimental setting, related terminology, and nomenclature.) The results confirm the conclusions drawn above for the case of the 30-bus model---operational variables should be taken into account because ignoring them leads to a significant increase in the investment cost, or even worse, infeasibility of a highly loaded configuration.

\begin{table}[h!]
\caption{Monetary advantage of considering operational variables for the 2736-bus Polish model.}
\centering
\begin{tabu} to \columnwidth { | X[c] | X[c] | X[c] |}
\hline
 Oper. variables  & Investment cost (\$) & Operational cost (\$/hour) \\
\hline
Fixed  & 916616.1 & 1884214.9  \\
\hline
Free   & 187869.6  & 1950027.2 \\
\hline
\end{tabu}
\label{tab5}
\end{table}

\paragraph{Importance of including operational cost in the objective function}

The normal operation base-case is taken with all loads re-scaled up by $5\%$, and the resulting optimizations are compared (including and not including the operational cost in the objective). The comparison is made for the total cost accumulated in 10 years. It is observed that the difference between the 0-year case (where the operational cost is ignored) and the 10-year case is 2.6$\%$, which results in 4330 M$\$$ of total cost savings; the additional investment (installation) cost is only 550k$\$$. The numbers clearly support main hypothesis: installation is advantageous and including the operation cost in the objective is mandatory for practical grid extension planning. This is possible because congestion in the system shows a decrease, in addition to the restoration of the feasible solution.

\paragraph{Scalability}

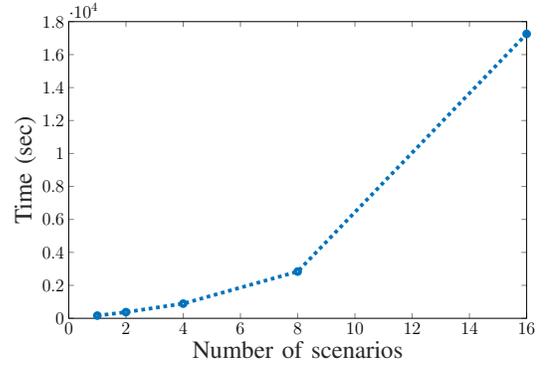
\begin{figure}[h!]
\begin{center}
\resizebox{0.8\columnwidth}{!}{
%
%
\definecolor{mycolor1}{rgb}{0.00000,0.44700,0.74100}%
\begin{tikzpicture}

\begin{axis}[%
font=\Large,
width=5.134in,
height=3.329in,
at={(0.861in,0.45in)},
scale only axis,
xmin=0,
xmax=16,
xlabel style={font=\color{white!15!black}},
xlabel={\huge Number of scenarios},
ymin=0,
ymax=18000,
ylabel style={font=\color{white!15!black}},
ylabel={\huge Time (sec)},
axis background/.style={fill=white}
]
\addplot [color=mycolor1, dashed, line width=3.0pt, mark=o, mark options={solid, mycolor1}, forget plot]
  table[row sep=crcr]{%
1	164.5\\
2	375.23\\
4	890.62\\
8	2837\\
16	17256\\
};
\end{axis}
\end{tikzpicture}%
}
\caption{Computational time vs the number of scenarios for the Polish model.}
\label{timee_pol}
\end{center}
\end{figure}

Fig.~\ref{timee_pol} shows how the computational time of the algorithm scales with the number of the scenarios.

A number of cases ranging from a single scenario to 16 scenarios are tested.  The brute-force IPOPT fails to solve the Polish model case even with a single scenario (exact problem in Eq.~\ref{OPT}). Proposed algorithm solves the most challenging case of 16 scenarios in 17500 sec. It is deduced from Fig.~\ref{timee_pol} that the computational time grows polynomially as {$O(K^3)$}, suggesting that proposed algorithm is practical/scalable for planning problem when computational time is not a significant constraint.

Note that the $O(K^3)$ scaling is still slower (with the $K$ increase) than the linear-scaling behavior observed in the 30-bus model. Our suggestion is that the better performance observed in the 30-bus model may be related to the fact that the Polish model is denser,  thus requiring linearization of more PF equations. It may also be due to worse scaling of the QP solver performance in the case of the Polish model. We plan to perform more detailed analysis in the future.

\section{Numerical Experiments: Multiple-Scenario-Aware Long-term Planning}
\label{sec:planning}

In this section developed methodology and algorithm, discussed in the preceding sections, are applied to analysis of the comprehensive multiple-scenario-aware long-term planning setting.  To generate the scenarios and initialize the algorithm, methods discussed in Appendix \ref{sec:scenarios} and Appendix \ref{sec:init_gen} are utilized. In all experiments discussed below the planning horizon is chosen to be 10 years.

\subsection{IEEE 30-bus model}
\begin{figure}[h]
\centering
\includegraphics[width=0.7\columnwidth]{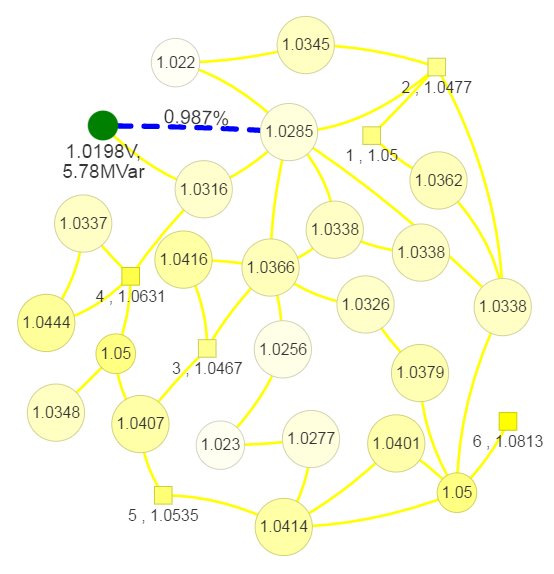}
\caption{Optimal solution for 10 years of planning in the case of the 30-bus model. Loads = yellow circles; gens = squares; blue dashed line = line with installed SC that was overloaded initially for some scenarios; green dot = node where an SVC is installed. Voltage levels are shown in PU; capacities of SVCs and SCs are shown in MVar and in $\%$ of initial line inductance.
\label{fig:10yearsLDcurveInvestCase30}}
\end{figure}

Sixteen scenarios per yearly LD curve are generated (160 total). The annual increase factor, $\beta$, is set to 1.5\% a year. The resulting optimal solution is shown in Fig.~\ref{fig:10yearsLDcurveInvestCase30}. It is observed that proposed algorithm installs FACTS devices efficiently and sparsely,  thus resolving successfully the otherwise imminent (observed for a significant portion of the 160 scenarios) AC-OPF infeasibility.   The optimal solution consists of the installation of an SVC device at bus 8 with the capacity of 5.78~MVAr and installation of an SC device at the line between buses 6 and 8 with a capacity increase of 1\%. The proposed investment is 30k$\$$, resulting in an average savings of 1.8\$ per hour.\footnote{All the actual values are model dependent. Costs are values of the objective function - determined according to a pre-defined cost of investement for a unit of capacity.}

\begin{figure} [h]
\centering
\includegraphics[width=0.9\columnwidth]{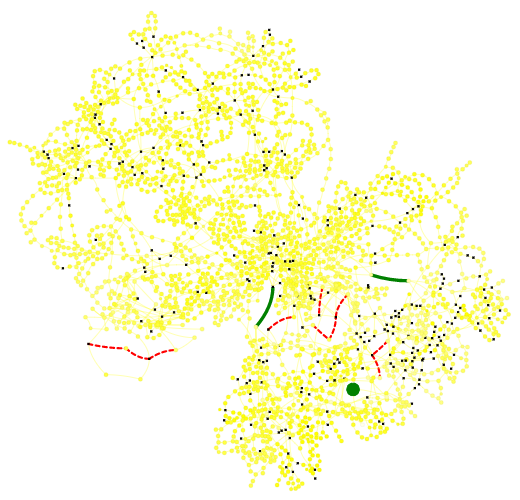}
\caption{Optimal solution for 10 years of planning in the case of the Polish model. Red dashed lines = lines that were initially overloaded for some scenarios. Two built SCs are shown by green lines. The big green dot illustrates an SVC device. SVC: 3.3 MVar, SCs: 14.4$\%$, 70.4$\%$ (left to right).
\label{fig:10yearsLDcurveInvestPol}}
\end{figure}

\subsection{2736-bus Polish model}
This experiment is done with 16 sampled scenarios (2 of 16 are AC-OPF infeasible) for the 10-year horizon and with assumed yearly loading growth (factor $\beta$) of 0.5\%. The resulting optimal investment is shown in Fig.~\ref{fig:10yearsLDcurveInvestPol} (coding of loads and gens is the same). The algorithm outputs a solution resulting in installation of two SC devices and one SVC device to resolve the infeasibility of some of the samples. Average congestion cost of the sampled scenarios is 5738\$/hour, and the average generation cost savings is 3369\$/hour. The solution is sparse and nonlocal (new FACTS devices are installed sufficiently far away from nodes and lines where the initial congestion was observed). SVC installation is relatively small in that case because congestion was much more significant for sampled scenarios than infeasibility, basically sampled two AC-OPF infeasible conditions are still close to the feasibility region.

\section{Conclusion}
\label{sec:conclusions}

New optimization framework for placing and sizing FACTS devices is proposed in this manuscript. The framework takes into account AC-PF equations. The most important features of the newly developed framework are the \textbf{scalability} of the algorithm, allowing to resolve congestion over practical (thousands of buses) size transmission systems, and the ability of the algorithm to \textbf{resolve multiple scenarios simultaneously}. Introduced optimization setting can also be considered as generalizing the standard AC-OPF: it seeks a balance between \textbf{installation} and \textbf{operations}.  The optimization objective includes the cost of operations over an extended time horizon as well as the cost of installation of the FACTS devices, represented in the form of an $l_1$ norm to promote \textbf{sparsity} of the resulting solutions. Optimization variables include capacities of the FACTS devices and respective operational settings associated with each scenario, where the latter are bounded by the former.

Proposed solution algorithm was tested in different regimes on a midsize model (30-bus IEEE) and a realistic size (Polish grid) model. It is observed that the output is spatially \textbf{sparse}, i.e., a very small number of FACTS devices is sufficient, and that the output is \textbf{nonlocal}, i.e., a typical new installation resolves congestion at multiple locations that can be rather far from the newly installed devices. Also it is observed that under highly loaded conditions FACTS devices are beneficial in reducing the total cost of generation. Optimal installation of the devices helps to resolve infeasibilities that are projected to become even more severe in the future as well.

The main technical achievement reported in this manuscript is the development of the algorithm that constitutes an efficient heuristic for solving the nonlinear and nonconvex optimization. The algorithm is sequential---it constructs a convergent sequence of convex and analytic formulations (QP with linear constraints) where  each constraint is represented explicitly through exact/analytic linearization of the original nonlinear constraints  (e.g., representing power balance at nodes and apparent power line limits) over all the degrees of freedom (including FACTS corrections) around the current operational point.

To represent uncertainty in the expected growth of the system (loads), a scenario sampling methodology is introduced. It is evident from experimental results that developed framework/approach is of a practical value for planning transmission grid expansion---it simultaneously resolves the growing economy and emerging congestion. Development of a convenient and flexible software for web visualization of transmission system states and FACTS installations has become a side benefit of this project.

Our current path forward is improvement of uncertainty modelling (generation side uncertainty related to renewables). For instance we are working on implementation of stochastic optimization methods and their applications to investment planning for the improvement of the value of the developed methodology and tools.

\appendices

\section{Transmission line modeling}
\label{sec:transmission_modeling}

For the sake of completeness and better understanding of the proposed heuristics, a short description of the  transmission line modeling is presented here. The so-called $\pi$-model from \cite{matpower} is utilized. Parameters of the line model are the series impedance, $z = r + \myspecial{j} x$, the total charging susceptance, $b$, the transformation ratio, $\tau$,  and the shift angle $\theta_{\mbox{shift}}$. A transformer breaks the symmetry between the \enquote{from} end, positioned next to the transformer, and the \enquote{to} end of the line.

Explicit expressions for apparent powers injected at the \enquote{from} end and the \enquote{to} end of the line in terms of voltages and phases are
\begin{eqnarray}
S_{f}(V_{f},\theta_{f},V_{t},\theta_{t},x)=\frac{V_{f}\left(r V_{f}- \tau V_{t}\left( r\cos\Delta+ x \sin\Delta\right)\right)}{\tau^2 l}\nonumber\\
-\myspecial{j} \frac{V_{f}}{2\tau^2 l}\Biggl(V_{f} (-2x+b l)+2\tau V_{t}\left(x \cos\Delta+r \sin\Delta\right))\Biggr) \label{Sf}
\end{eqnarray}
\begin{eqnarray}
S_{t}(V_{f},\theta_{f},V_{t},\theta_{t},x)=\frac{V_{t} \left(r \tau V_{t}-V_{f} (r \cos\Delta+x \sin\Delta)\right)}{\tau^2 l}\nonumber \\
-\myspecial{j}\frac{V_{t}}{2\tau l}\Biggl(\tau V_{t} (-2x+bl)+2V_{f} (x\cos\Delta-r\sin\Delta)\Biggr)\label{St}
\end{eqnarray}
where $\forall i:\quad v_i=V_i e^{{\it \bf j}\theta_i}$, $\Delta = \theta_{f}-\theta_{t}-\theta_{\mbox{shift}}$ and $l=r^2+x^2$.

\section{Generation of scenarios}
\label{sec:scenarios}

\begin{figure}[h!]
\centering
\def\svgwidth{0.6\columnwidth}
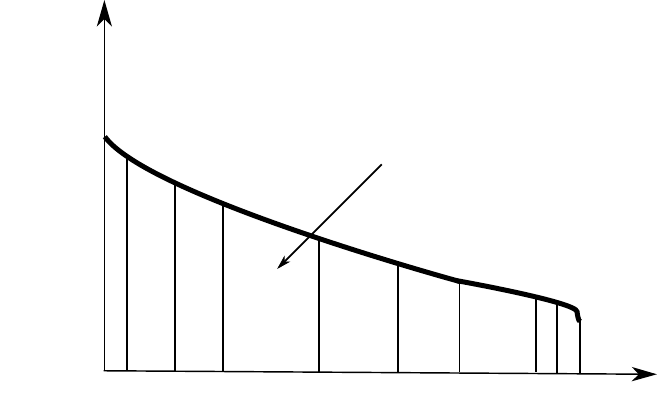
\caption{Piecewise-constant approximation of the LD curve.
\label{fig:LDcurveAppr}}
\end{figure}

Scenario generation/sampling is used to include the uncertainty related to system load for the planning period. Power system load growth over the time horizon is modeled via modification of the Load Duration (LD) curve for current year. The base LD curve is illustrated in Fig.~\ref{fig:LDcurveAppr}.

The base LD curve is used first to generate LD curves for consecutive years, rescaling the base LD curve by the load growth factor of $0.5\%-1.5\%$ a year. Second, each early LD curve is split into $M$ piecewise-constant parts. ($M=6$ in simulations.) Finally,  each piece of an LD curve is used to generate scenarios according to a random (thus called sampling) procedure described below. This scheme of scenario generation/sampling models variations in the distribution of loads, thus  simulating power system behavior during an extended period of time in the future.

It is assumed (and this assumption is confirmed in all our experimental tests) that each of the generated (sampled) load scenarios is ACOPF-feasible when the line constraints are ignored. (In other words, the setting is considered when there is enough generation capacity even for the stressed cases.) Depending on the sampled scenario, three situations may arise.
\begin{enumerate}
\item ACOPF is feasible and congestion price is zero (low loading level).
\item ACOPF is feasible and congestion price is positive (higher loading level representing peak conditions).
\item ACOPF is infeasible because of congestion of lines and/or voltage constraints but the system has enough generation capacity. ACOPF without apparent power limits on lines (and without voltage constraints if infeasible) is feasible (overloaded conditions that are possible in the future).
\end{enumerate}
The aim of planning installation of FACTS devices at the right locations with their corresponding capacities is to reduce generation cost for situation 2 and to improve or extend the feasibility domain of the system for situation 3. Extra years of service can hence be added to the existing grid by making it more flexible, thereby delaying investments into new lines and generators.

\subsection{Scenario sampling for each segment}
\label{subsec:scenarios_sampling}

The loading level $\alpha_i$ for a segment $i$ is represented by
\begin{equation}
\alpha_i=\frac{\overline{\alpha}_i+\underline{\alpha}_i}{2}
\end{equation}
Future loading configurations are obtained from the base case by rescaling all active and reactive loads by $\alpha_i$ uniformly. The resulting vector of loads for a segment is thus given by
\begin{eqnarray}
    l_i^0=\alpha_i \times l^0
\end{eqnarray}
Loading configurations are generated for each segment $i$ and each $j=1..N_i$ through modification of initial $l_i^0$. This is done by adding Gaussian correction to each load with zero expected value and a respective standard deviation:
\begin{eqnarray}
   & & l_i^j = l_i^0 + \mathcal{N}(0, \sigma_{l_i^0})  \label{Gauss}\\
   & & p_i^j = w_i/N_i \quad (\mbox{probability of a given scenario}) \label{Gauss-p}
\end{eqnarray}
where $\sigma_{l_i^0}$ is given by
 \begin{align}
  \sigma_{l_i^0}&=\frac{\overline{\alpha}_i-\alpha_i}{\alpha_i} \times l_i^0	\\
	& =\sigma \times l_i^0
 \end{align}
The choice of parameters used in our experimental test to sample the scenarios is described in Table~\ref{fig:ActualLDappr}.

\begin{table} [h!]
\caption{Implementation of the LD curve scheme}
\centering
\begin{tabu} to \columnwidth { | X[c] | X[c] | X[c] | X[c] |}
\hline
$i$ & $w_i$ & $\alpha_i$ & $\sigma$ \\ [0.5ex]
\hline
1 & 5.50 & 0.940 & 0.064 \\
2 & 19.50 & 0.845 & 0.041 \\
3 & 25.00 & 0.775 & 0.045 \\
4 & 25.00 & 0.685 & 0.080 \\
5 & 18.80 & 0.590 & 0.068 \\
6 & 6.20 & 0.51 & 0.078 \\
\hline
\end{tabu}
\label{fig:ActualLDappr}
\end{table}

\subsection{Congestion analysis correction}
\label{subsec:congestion}

If the case is considered in which, for a given load configuration, standard ACOPF outputs a solution that is not congested,  i.e., a solution for which each  constraint (on line flows or voltages) is satisfied with a  margin,  then this scenario does not require any FACTS device installation. If the whole segment (from the procedure described in the preceding subsection) is of this ``zero-congestion" type,  then obviously one does not need to generate many samples for the segment.  Instead, one rescaled base scenario to represent the whole segment is picked.

\section{Defining initial generation profile}	
\label{sec:init_gen}

The initial profile of the generation for each load scenario has to be determined to run the algorithm. Generation capacity is assumed to be large enough for given loading levels. Two procedures are used for that: (1) solve ACOPF with the thermal limits ignored, and (2) find proportional generator response. Second is done in the following four steps:
\begin{itemize}
\item Search for the smallest load rescaling factor $\alpha$ lowering the load and thus making the resulting case feasible.
\item Solve ACOPF with this new rescaled loading.
\item Proportionally increase generation and load with the value of $\alpha$, which restores the initial loading of the system. Use voltages from the ACOPF solution.
\item Solve ACPF to obtain generation maintaining the loading.
\end{itemize}

\ifCLASSOPTIONcaptionsoff
  \newpage
\fi

\bibliography{bibl_new}

\vskip 0pt plus -1fil

\begin{IEEEbiography}
[{\includegraphics[width=1in,height=1.25in,clip,keepaspectratio]{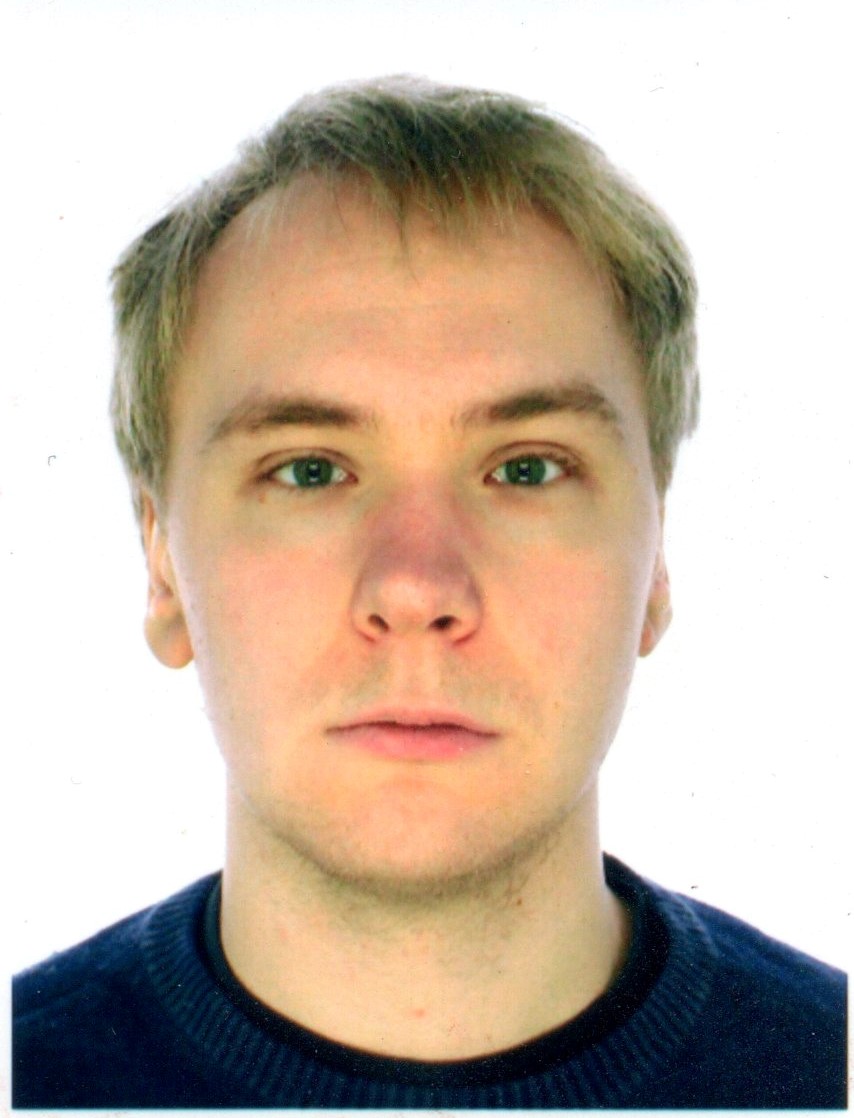}}]
{Vladimir Frolov} was born in Vologda, Russia in 1991. He received the B.S. and M.S. degrees in applied physics and mathematics from Moscow Institute of Physics and Technology, Moscow, in 2015. He also received the M.S. degree in energy science and technology from Skolkovo Institute of Science and Technology (Skoltech), Moscow, in 2015 and he is currently on his final year of the Ph.D. degree in engineering systems at Skoltech.
His current research interests include power system modelling, nonlinear optimization, operational and investment planning under uncertainty and power system reinforcement.
\end{IEEEbiography}

\vskip 0pt plus -1fil

\begin{IEEEbiography}
[{\includegraphics[width=1in,height=1.25in,clip,keepaspectratio]{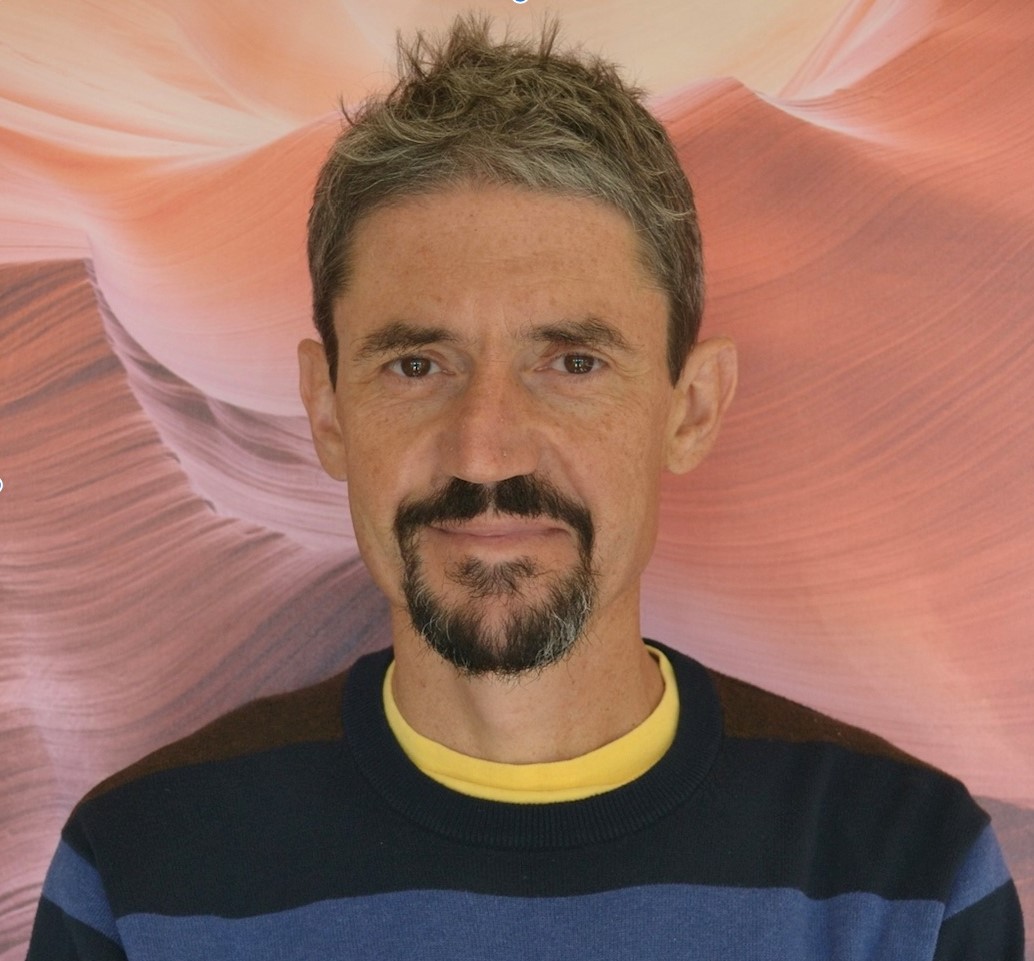}}]
{Michael Chertkov}'s areas of interest and expertise include mathematics and statistics applied to physical, engineering and data sciences. He received his Ph.D. in physics from the Weizmann Institute of Science in 1996, and his M.Sc. in physics from Novosibirsk State University in 1990.

After his Ph.D., he spent three years at Princeton University as a R.H. Dicke Fellow in the Department of Physics. He joined Los Alamos National Lab in 1999, initially as a J.R. Oppenheimer Fellow in the Theoretical Division, and continued as a Technical Staff Member leading projects in physics of algorithms, energy grid systems, physics and engineering informed data science and machine learning for turbulence.

In 2019 Dr. Chertkov moved to Tucson to become professor of mathematics and lead Interdisciplinary Graduate Program in Applied Mathematics at the University of Arizona, continuing to work for LANL part time.  Dr. Chertkov has published more than 200 papers. He is a fellow of the American Physical Society (APS) and a senior member of IEEE.
\end{IEEEbiography}


\end{document}